# ON SLANT HELICES AND GENERAL HELICES IN EUCLIDEAN $n$-SPACE

Yusuf YAYLI[1], Evren ZIPLAR[2]

[1] *Department of Mathematics, Faculty of Science, University of Ankara, Tandoğan, Ankara, Turkey*

*yayli@science.ankara.edu.tr*

[2] *Department of Mathematics, Faculty of Science, University of Ankara, Tandoğan, Ankara, Turkey*

*evrenziplar@yahoo.com*

**Abstract.** In this paper, in Euclidean $n$-space $E^n$, we investigate the relation between slant helices and spherical helices. Moreover, in $E^n$, we show that a slant helix and the tangent indicatrix of the slant helix have the same axis (or direction). Also, we give the important relations between slant helices, spherical helices in $E^n$ and geodesic curves on a helix hypersurface in $E^n$.

*Keywords*: Slant helices; general helices; spherical helices; tangent indicatrix; helix hypersurfaces

**Mathematics Subject Classification 2000:** 53A04, 53B25, 53C40, 53C50

## 1. INTRODUCTION

Slant helix is one of the most important topics of differential geometry. Izumiya and Takeuchi have investigated the many properties of slant helices that the normal lines make a constant angle with a fixed direction in Euclidean 3-space [7]. Moreover, they proved that a space curve is a slant helix if and only if the geodesic curvature of the principal normal of the curve is a constant function [7]. Besides, Izumiya and Takeuchi described a slant helix if and only if the function

$$\frac{\kappa^2}{(\kappa^2+\tau^2)^{3/2}}(\frac{\tau}{\kappa})'$$

is constant, where $\kappa$ is curvature and $\tau$ is torsion of the curve, respectively[7].

Monterde obtains a geometric characterization of Salkowski curves whose the normal vector maintains a constant angle with a fixed direction in space [11].

On the other hand, Kula and Yayli have investigated spherical images of tangent indicatrix of a slant helix in Euclidean 3-space and they proved that the spherical images are spherical helix in Euclidean 3-space[9].

General helice whose tangents make a constant angle with a fixed direction is also considerable subject of differential geometry. Many geometers have studied on this type of curves [5,10,14]. In 1845, de Saint Venant first proved that a space curve is a general helix if and only if the ratio of curvature to torsion be constant[15]. Moreover, Ali and López consider the generalization of the concept of general helices in the Euclidean $n$-space $E^n$ [1]. Also, Yılmaz and Turgut introduce a new version of Bishop frame and introduce new spherical images [17].

In differential geometry of surfaces, an helix hypersurface in $E^n$ is defined by the property that tangent planes make a constant angle with a fixed direction [4]. Di Scala and Ruiz-Hernández have introduced the concept of these surfaces in [4]. And, A.I.Nistor has also introduced certain constant angle surfaces constructed on curves in $E^3$ [12].

Özkaldı and Yayli give some characterization for a curve lying on a surface for which the unit normal makes a constant angle with a fixed direction [13].

One of the main purpose of this work is to observe the relations between slant helices and spherical helices. Another purpose of this study is to give an important relation between helix hypersurfaces and spherical helices.

## 2. PRELIMINARIES

**Definition 2.1** *Let $\alpha : I \subset \mathrm{IR} \to E^n$ be an arbitrary curve in $E^n$. Recall that the curve $\alpha$ is said to be of unit speed ( or parametrized by the arc-length function $s$ ) if $\langle \alpha'(s), \alpha'(s) \rangle = 1$, where $\langle , \rangle$ is the standart scalar product in the Euclidean space $E^n$ given by*

$$\langle X, Y \rangle = \sum_{i=1}^{n} x_i y_i ,$$

*for each $X = (x_1, x_2, ..., x_n), Y = (y_1, y_2, ..., y_n) \in E^n$.*

Let $\{V_1(s), V_2(s),...,V_n(s)\}$ be the moving frame along $\alpha$, where the vectors $V_i$ are mutually orthogonal vectors satisfying $\langle V_i, V_i \rangle = 1$. The Frenet equations for $\alpha$ are given by

$$\begin{bmatrix} V_1' \\ V_2' \\ V_3' \\ \vdots \\ V_{n-1}' \\ V_n' \end{bmatrix} = \begin{bmatrix} 0 & k_1 & 0 & 0 & \cdots & 0 & 0 \\ -k_1 & 0 & k_2 & 0 & \cdots & 0 & 0 \\ 0 & -k_2 & 0 & k_3 & \cdots & 0 & 0 \\ \vdots & \vdots & \vdots & \vdots & & \vdots & \vdots \\ 0 & 0 & 0 & 0 & \cdots & 0 & k_{n-1} \\ 0 & 0 & 0 & 0 & \cdots & -k_{n-1} & 0 \end{bmatrix} \begin{bmatrix} V_1 \\ V_2 \\ V_3 \\ \vdots \\ V_{n-1} \\ V_n \end{bmatrix}.$$

Recall that the functions $k_i(s)$ are called the $i$-th curvatures of $\alpha$ [6].

**Definition 2.2** A unit speed curve $\beta: I \subset \mathrm{IR} \to \mathrm{E}^n$ is called general helix if its tangent vector $V_1^*$ makes a constant angle with a fixed direction $U$ [1].

**Theorem 2.1** Let $\beta: I \to \mathrm{E}^n$ be a unit speed curve and a general helix in $\mathrm{E}^n$. Then the unit direction of the helix :

$$U = \cos\theta_1 [V_1^* + \sum_{i=3}^{n} G_i^* V_i^*].$$

Here, $\langle V_1^*(s), U \rangle = \cos\theta_1$ and $G_1^* = 1, G_2^* = 0, G_i^* = \dfrac{1}{k_{i-1}^*}[k_{i-2}^* G_{i-2}^* + (G_{i-1}^*)']$, $3 \leq i \leq n$. On the other hand, $\{V_1^*, V_2^*,...,V_n^*\}$ is the Frenet frame of $\beta$ and $k_{i-1}^*, k_{i-2}^*$ are curvatures of $\beta$, where $3 \leq i \leq n$ [1].

**Definition 2.3** A unit speed curve $\alpha: I \subset \mathrm{IR} \to \mathrm{E}^n$ is called slant helix if its unit principal normal $V_2$ makes a constant angle with a fixed direction $L$ [2].

**Theorem 2.2** Let $\alpha: I \to \mathrm{E}^n$ be a unit speed curve in $\mathrm{E}^n$. Define the functions

$$G_1 = \int k_1(s) ds, \; G_2 = 1, \; G_3 = \dfrac{k_1}{k_2} G_1, \; G_i = \dfrac{1}{k_{i-1}}[k_{i-2} G_{i-2} + G_{i-1}'],$$

where $4 \leq i \leq n$. Then $\alpha$ is a slant helix if and only if the function

$$\sum_{i=1}^{n} G_i^2 = C$$

is constant and non-zero. Moreover, the constant $C = \sec^2 \theta_2$, being $\theta_2$ the angle that makes $V_2$ with the fixed direction $L$ that determines $\alpha$ [2].

**Theorem 2.3** *Let $\alpha : I \to E^n$ be a unit speed curve and a slant helix in $E^n$. Then the unit direction of the helix :*

$$L = \cos \theta_2 \, [\sum_{i=1}^{n} G_i V_i].$$

*Here, $\langle V_2(s), L \rangle = \cos \theta_2$ and*

$$G_1 = \int k_1(s)ds \, , G_2 = 1, G_3 = \frac{k_1}{k_2} G_1 \, , G_i = \frac{1}{k_{i-1}}[k_{i-2}G_{i-2} + G'_{i-1}], 4 \leq i \leq n.$$

*On the other hand, $\{V_1, V_2, ..., V_n\}$ is the Frenet frame of $\alpha$ [2].*

**Definition 2.4** *Let $\beta : I \to E^n$ be a unit speed curve in $E^n$. Harmonic curvatures of $\beta$ is defined by*

$$H_i : I \subset \mathbb{R} \to \mathbb{R} \, , i = 0,1,2,...,n-2,$$

$$H_i = \begin{cases} 0, & i = 0 \\ \dfrac{k_1}{k_2}, & i = 1 \\ \{V_1[H_{i-1}] + H_{i-2} k_i\}\dfrac{1}{k_{i+1}}, & i = 2,3,...,n-2 \end{cases}$$

*[3].*

**Theorem 2.4** *Let $\beta : I \to E^n$ be a unit speed curve and a general helix in $E^n$. Then the unit direction of the helix :*

$$X = \cos \theta_1 \, [V_1^* + H_1 V_3^* + ... + H_{n-2} V_n^*].$$

*Here, $\langle V_1^*, X \rangle = \cos \theta_1$. On the other hand, $\{V_1^*, V_2^*, ..., V_n^*\}$ is the Frenet frame of $\beta$ and $\{H_1, H_2, ..., H_{n-2}\}$ are the harmonic curvatures of $\beta$ [3].*

**Remark 2.1** If the Frenet frame of the tangent indicatrix $\beta$ of a space curve $\alpha$ is $\{V_1^* = T, V_2^* = N, V_3^* = B\}$, then

$$V_1^* = T = n$$

$$V_2^* = N = \frac{1}{\sqrt{\kappa^2 + \tau^2}}(-\kappa t + \tau b)$$

$$V_3^* = B = \frac{1}{\sqrt{\kappa^2 + \tau^2}}(\tau t + \kappa b)$$

and $\kappa_\beta = \frac{\sqrt{\kappa^2 + \tau^2}}{\kappa}$ is the curvature of $\beta$, $\tau_\beta = \frac{\kappa \tau' - \kappa' \tau}{\kappa(\kappa^2 + \tau^2)}$ is the torsion of $\beta$, where $\{V_1 = t, V_2 = n, V_3 = b\}$ is the Frenet frame of $\alpha$ and $\kappa$ is the curvature of $\alpha$, $\tau$ is the torsion of $\alpha$ [8].

## 3. SLANT HELICES AND SPHERICAL HELICES

In the following Theorem, we give the relation between slant helices in $E^n$ and spherical helices on the unit hypersphere in $S^{n-1} \subset E^n$.

**Theorem 3.1** Let $\alpha : I \subset \mathrm{IR} \to E^n$ be a unit speed curve (parametrized by arclength function $s$) in $E^n$ and let $\beta : I \to S^{n-1} \subset E^n$ be the tangent indicatrix of the curve $\alpha$, where $S^{n-1}$ is the unit hypersphere in $E^n$. Then the curve $\alpha$ is a slant helix with direction $L$ in $E^n$ if and only if the curve $\beta$ is a general helix (spherical helix) with direction $L$ on $S^{n-1} \subset E^n$. In other words, $\alpha$ and $\beta$ have the same direction $L$.

**Proof:** The tangent indicatrix of the curve $\alpha$ is defined by

$$\frac{d\alpha}{ds} = \beta(s).$$

We assume that the arclength parameter of $\beta$ is $s_\beta$. Then we can write

$$\frac{d\beta}{ds_\beta} = \frac{d\beta}{ds} \frac{ds}{ds_\beta}.$$

And, from the Frenet equations (see Definition 2.1),

$$\frac{d\beta}{ds_\beta} = k_1(s)V_2(s)\frac{ds}{ds_\beta}.$$

Thus, from the last equation, by taking norms on both sides, we obtain $\frac{ds}{ds_\beta} = \frac{1}{k_1(s)}$

$(k_1(s) \neq 0)$. Hence, we have $\frac{d\beta}{ds_\beta} = V_2(s)$.

Now, we assume that the curve $\alpha$ is a slant helix with direction $L$ in $E^n$. So, the scalar product between the unit principal normal vector field $V_2$ with $L$ is

$$\langle V_2, L \rangle = \cos\theta.$$

On the other hand, we know that $\frac{d\beta}{ds_\beta} = V_2$. Hence, we have $<\frac{d\beta}{ds_\beta}, L> = \cos\theta$. It follows that the curve $\beta$ is a general helix (spherical helix) with the direction $L$ on $S^{n-1} \subset E^n$.

Conversely, the curve $\beta$ is a general helix (spherical helix) with direction $L$ on $S^{n-1} \subset E^n$. Then we have $<\frac{d\beta}{ds_\beta}, L> = \cos\theta$, where $s_\beta$ is the arclength parameter of $\beta$.

On the other hand, we know that $V_2 = \frac{d\beta}{ds_\beta}$. Hence, we have $\langle V_2, L \rangle = \cos\theta$. So, we deduce that the curve $\alpha$ is a slant helix with direction $L$ in $E^n$. This completes the proof.

This above Theorem has the following corollary.

**Corollary 3.1** *Let $\alpha: I \subset \mathrm{IR} \to E^n$ be a slant helix with unit speed (or parametrized by arclength function $s$) and let $\beta: I \to S^{n-1} \subset E^n$ be the tangent indicatrix of the curve $\alpha$ (parametrized by arclength function $s_\beta$), where $S^{n-1}$ is the unit hypersphere in $E^n$. We assume that the direction of $\alpha$ is $L = \cos\theta[\sum_{i=1}^{n} G_i V_i]$. Then the direction $L$ can be expressed in the forms*

$$L = \cos\theta [V_1^* + \sum_{i=3}^{n} G_i^* V_i^*] \text{ and } L = \cos\theta [V_1^* + H_1 V_3^* + \ldots + H_{n-2} V_n^*],$$

*where $\theta$ is constant.*

*Here, for the curve $\alpha$:*

$$G_1 = \int k(s)ds, \ G_2 = 1, \ G_3 = \frac{k_1}{k_2} G_1, \ G_i = \frac{1}{k_{i-1}}[k_{i-2} G_{i-2} + G'_{i-1}], \ 4 \leq i \leq n$$

*and $\{V_1, V_2, \ldots, V_n\}$ is the Frenet frame of $\alpha$.*

*Moreover, for the curve $\beta$:*

$$G_1^* = 1, \ G_2^* = 0, \ G_i^* = \frac{1}{k_{i-1}^*}[k_{i-2}^* G_{i-2}^* + (G_{i-1}^*)'], \ 3 \leq i \leq n$$

*and $\{V_1^*, V_2^*, \ldots, V_n^*\}$ is the Frenet frame of $\beta$, $\{H_1, H_2, \ldots, H_{n-2}\}$ are the harmonic curvatures of $\beta$.*

**Proof:** From the Theorem 3.1, since $\alpha$ and $\beta$ have the same direction $L$ and $\frac{d\beta}{ds_\beta} = V_1^* = V_2$, then $\langle V_1^*, L \rangle = \langle V_2, L \rangle = \cos\theta$. So, we deduce that

$$\cos\theta [\sum_{i=1}^{n} G_i V_i] = \cos\theta [V_1^* + \sum_{i=3}^{n} G_i^* V_i^*] = \cos\theta [V_1^* + H_1 V_3^* + \ldots + H_{n-2} V_n^*].$$

This completes the proof.

In the following example, it has been obtained a general helix on the unit hypersphere $S^2 \subset E^3$ by using a slant helix in $E^3$.

**Example 3.1** Let $\alpha(s) = (\frac{2}{5}\sin 2s - \frac{1}{40}\sin 8s, \frac{-2}{5}\cos 2s + \frac{1}{40}\cos 8s, \frac{4}{15}\sin 3s)$ be a slant helix with unit speed in $E^3$, where $\pi/3 < s < 2\pi/3$. It is easily obtain the curvatures as follows:

$$\kappa(s) = -4\sin 3s \text{ and } \tau(s) = 4\cos 3s,$$

where, $\kappa$ is the curvature of $\alpha$ and $\tau$ is the torsion of $\alpha$.

Now, we are going to find out the tangent indicatrix $\beta$ of the curve $\alpha$:

$$\beta(s) = \frac{d\alpha}{ds} = (\frac{4}{5}\cos 2s - \frac{1}{5}\cos 8s, \frac{4}{5}\sin 2s - \frac{1}{5}\sin 8s, \frac{4}{5}\cos 3s).$$

That is, $\beta(s)$ is a general helix on the hypersphere $S^2$ in $E^3$. The slant helix $\alpha$ is shown in the following Figure 1 and the tangent indicatrix $\beta$ of the curve $\alpha$ is shown in the following Figure 2, where $\pi/3 < s < 2\pi/3$.

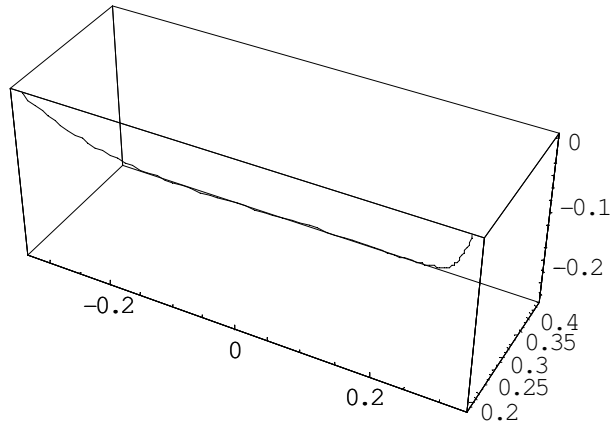

Figure 1

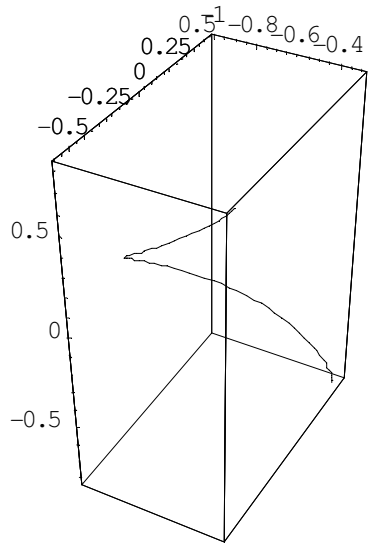

Figure 2

For Corollary 3.1, we can give the following example.

**Example 3.2** In this example, we are going to find out the directions of the curves $\alpha$ and $\beta$ defined in Example 3.1.

The direction of the slant helix $\alpha$:

$$L = \cos\theta [\sum_{i=1}^{3} G_i V_i] = \cos\theta [G_1 t + G_2 n + G_3 b],$$

where $G_1 = \int \kappa(s) ds$, $G_2 = 1$, $G_3 = \dfrac{\kappa}{\tau} G_1$ and $\{V_1 = t, V_2 = n, V_3 = b\}$ is the Frenet frame of the slant helix $\alpha$. From the example 3.1, we know that $\kappa(s) = -4\sin 3s$ and $\tau(s) = 4\cos 3s$. So, it is easily obtain the $G_1$ and $G_3$ as follows:

$$G_1 = \frac{4}{3}\cos 3s \text{ and } G_3 = \frac{-4}{3}\sin 3s.$$

On the other hand, from theorem 2.2, we know that $G_1^2 + G_2^2 + G_3^2 = \sec^2\theta$. Hence, we can deduce from the last equality that $\cos\theta = \dfrac{3}{5}$.

Finally, the direciton of the slant helix $\alpha$ is found out as

$$L = (\frac{4}{5}\cos 3s)\, t + \frac{3}{5}n + (\frac{-4}{5}\sin 3s) b.$$

The direction of the tangent indicatrix $\beta$:

$$U = \cos\theta [V_1^* + \sum_{i=3}^{3} G_i^* V_i^*] = \cos\theta [T + G_3^* B],$$

where $G_3^* = \dfrac{1}{\tau_\beta}[\kappa_\beta G_1^* + (G_2^*)']$ ($\kappa_\beta$ is the curvature of $\beta$ and $\tau_\beta$ is the torsion of $\beta$) and $\{V_1^* = T, V_2^* = N, V_3^* = B\}$ is the Frenet frame of $\beta$.

From theorem 2.1, we know that $G_1^* = 1$ and $G_2^* = 0$. So, $G_3^* = \dfrac{\kappa_\beta}{\tau_\beta}$. On the other hand, from remark 2.1, $\kappa_\beta = \dfrac{\sqrt{\kappa^2 + \tau^2}}{\kappa}$ and $\tau_\beta = \dfrac{\kappa\tau' - \kappa'\tau}{\kappa(\kappa^2 + \tau^2)}$. Hence, we obtain

$G_3^* = \dfrac{(\kappa^2 + \tau^2)^{3/2}}{\kappa\tau' - \kappa'\tau}$. Again, from remark 2.1, $T = n$ and $B = \dfrac{1}{\sqrt{\kappa^2 + \tau^2}}(\tau t + \kappa b)$.

Consequently, if we also consider that $\kappa(s) = -4\sin 3s$ and $\tau(s) = 4\cos 3s$, then the direciton of the helix $\beta$ is found out as

$$U = (\frac{4}{5}\cos 3s)t + \frac{3}{5}n + (\frac{-4}{5}\sin 3s)b,$$

where $\cos\theta = \frac{3}{5}$.

Notice that the directions of the curves $\alpha$ and $\beta$ equal.

## 4. HELIX HYPERSURFACES AND HELICES

In this section, we give the important relations between slant helices, spherical helices in $\mathrm{E}^n$ and geodesic curves on a helix hypersurface in $\mathrm{E}^n$.

**Definition 4.1** *Given a hypersurface $M \subset \mathrm{E}^n$ and an unitary vector $d \neq 0$ in $\mathrm{E}^n$, we say that $M$ is a helix hypersurface with respect to the fixed direction $d$ if $\langle d, \xi \rangle$ is constant function along $M$, where $\xi$ is a normal vector field on $M$ [4].*

**Theorem 4.1** *Let $M$ be a helix hypersurface with the direction $d$ in $\mathrm{E}^n$ and let $\alpha : \mathrm{I} \subset \mathrm{IR} \to M$ be a unit speed geodesic curve on $M$. Then, the curve $\alpha$ is a slant helix with the direction $d$ in $\mathrm{E}^n$ [16].*

**Proof:** Let $\xi$ be a normal vector field on $M$. Since $M$ is a helix hypersurface with respect to $d$, $\langle d, \xi \rangle = $ constant. That is, the angle between $d$ and $\xi$ is constant on every point of the surface $M$. And, $\alpha''(s) = \lambda \xi \mid_{\alpha(s)}$ along the curve $\alpha$ since $\alpha$ is a geodesic curve on $M$. Moreover, by using the Frenet equation $\alpha''(s) = V_1' = k_1 V_2$, we obtain $\lambda \xi \mid_{\alpha(s)} = k_1 V_2$, where $k_1$ is first curvature of $\alpha$. Thus, from the last equation, by taking norms on both sides, we obtain $\xi = V_2$ or $\xi = -V_2$. So, $\langle d, V_2 \rangle$ is constant along the curve $\alpha$ since $\langle d, \xi \rangle = $ constant. In other words, the angle between $d$ and $V_2$ is constant along the curve $\alpha$. Consequently, the curve $\alpha$ is a slant helix with the direction $d$ in $E^n$.

This completes the proof.

Now, we can give a new Theorem by using the above Theorem.

**Theorem 4.2** *Let $\alpha: I \subset \mathbb{R} \to M \subset E^n$ be a geodesic curve on $M$ with unit speed (or parametrized by arclength function $s$), where $M$ is a helix hypersurface with respect to a fixed direction $d$ in $E^n$. Then the tangent indicatrix $\alpha'(s)$ of the curve $\alpha(s)$ is a spherical helix on unit the hypersphere $S^{n-1} \subset E^n$.*

**Proof:** We assume that $\alpha$ is a geodesic curve on $M$. Then from Theorem 4.1, $\alpha$ is a slant helix with the direction $d$ in $E^n$. On the other hand, from Theorem 3.1, the tangent indicatrix of the curve $\alpha$ is a spherical helix with the direction $d$ on the hypersphere $S^{n-1} \subset E^n$.

This completes the proof.

This above Theorem has the following corollary.

**Corollary 4.1** *Let $M \subset E^n$ be a helix hypersurface in $E^n$. Then, the tangent indicators of all geodesic curves on the surface $M$ are spherical helices whose axes coincide.*